\documentclass{amsart}
\title{An ``Anti-Hasse Principle'' for Prime Twists}
\author{Pete L. Clark}
\email{pete@math.uga.edu}

\newtheorem{lemma}{Lemma}

\newtheorem{cor}[lemma]{Corollary}
\newtheorem{thm}[lemma]{Theorem}

\DeclareFontEncoding{OT2}{}{} 
  \newcommand{\textcyr}[1]{%
    {\fontencoding{OT2}\fontfamily{wncyr}\fontseries{m}\fontshape{n}%
     \selectfont #1}}
\newcommand{\Sha}{{\mbox{\textcyr{Sh}}}}

\begin{document}
\maketitle
\newcommand{\DD}{\mathcal{D}}
\newcommand{\F}{\ensuremath{\mathbb F}}
\newcommand{\Fp}{\ensuremath{\F_p}}
\newcommand{\Fl}{\ensuremath{\F_l}}
\newcommand{\Fpbar}{\overline{\Fp}}
\newcommand{\Fq}{\ensuremath{\F_q}}
\newcommand{\PP}{\mathbb{P}}
\newcommand{\PPone}{\mathfrak{p}_1}
\newcommand{\PPtwo}{\mathfrak{p}_2}
\newcommand{\PPonebar}{\overline{\PPone}}
\newcommand{\N}{\ensuremath{\mathbb N}}
\newcommand{\Q}{\ensuremath{\mathbb Q}}
\newcommand{\Qbar}{\overline{\Q}}
\newcommand{\R}{\ensuremath{\mathbb R}}
\newcommand{\Z}{\ensuremath{\mathbb Z}}
\newcommand{\SSS}{\ensuremath{\mathcal{S}}}
\newcommand{\Rn}{\ensuremath{\mathbb R^n}}
\newcommand{\Ri}{\ensuremath{\R^\infty}}
\newcommand{\C}{\ensuremath{\mathbb C}}
\newcommand{\Cn}{\ensuremath{\mathbb C^n}}
\newcommand{\Ci}{\ensuremath{\C^\infty}}
\newcommand{\U}{\ensuremath{\mathcal U}}
\newcommand{\gn}{\ensuremath{\gamma^n}}
\newcommand{\ra}{\ensuremath{\rightarrow}}
\newcommand{\fhat}{\ensuremath{\hat{f}}}
\newcommand{\ghat}{\ensuremath{\hat{g}}}
\newcommand{\hhat}{\ensuremath{\hat{h}}}
\newcommand{\covui}{\ensuremath{\{U_i\}}}
\newcommand{\covvi}{\ensuremath{\{V_i\}}}
\newcommand{\covwi}{\ensuremath{\{W_i\}}}
\newcommand{\Gt}{\ensuremath{\tilde{G}}}
\newcommand{\gt}{\ensuremath{\tilde{\gamma}}}
\newcommand{\Gtn}{\ensuremath{\tilde{G_n}}}
\newcommand{\gtn}{\ensuremath{\tilde{\gamma_n}}}
\newcommand{\gnt}{\ensuremath{\gtn}}
\newcommand{\Gnt}{\ensuremath{\Gtn}}
\newcommand{\Cpi}{\ensuremath{\C P^\infty}}
\newcommand{\Cpn}{\ensuremath{\C P^n}}
\newcommand{\lla}{\ensuremath{\longleftarrow}}
\newcommand{\lra}{\ensuremath{\longrightarrow}}
\newcommand{\Rno}{\ensuremath{\Rn_0}}
\newcommand{\dlra}{\ensuremath{\stackrel{\delta}{\lra}}}
\newcommand{\pii}{\ensuremath{\pi^{-1}}}
\newcommand{\la}{\ensuremath{\leftarrow}}
\newcommand{\gonem}{\ensuremath{\gamma_1^m}}
\newcommand{\gtwon}{\ensuremath{\gamma_2^n}}
\newcommand{\omegabar}{\ensuremath{\overline{\omega}}}
\newcommand{\dlim}{\underset{\lra}{\lim}}
\newcommand{\ilim}{\operatorname{\underset{\lla}{\lim}}}
\newcommand{\Hom}{\operatorname{Hom}}
\newcommand{\Ext}{\operatorname{Ext}}
\newcommand{\Part}{\operatorname{Part}}
\newcommand{\Ker}{\operatorname{Ker}}
\newcommand{\im}{\operatorname{im}}
\newcommand{\ord}{\operatorname{ord}}
\newcommand{\unr}{\operatorname{unr}}
\newcommand{\B}{\ensuremath{\mathcal B}}
\newcommand{\Ocr}{\ensuremath{\Omega_*}}
\newcommand{\Rcr}{\ensuremath{\Ocr \otimes \Q}}
\newcommand{\Cptwok}{\ensuremath{\C P^{2k}}}
\newcommand{\CC}{\ensuremath{\mathcal C}}
\newcommand{\gtkp}{\ensuremath{\tilde{\gamma^k_p}}}
\newcommand{\gtkn}{\ensuremath{\tilde{\gamma^k_m}}}
\newcommand{\QQ}{\ensuremath{\mathcal Q}}
\newcommand{\I}{\ensuremath{\mathcal I}}
\newcommand{\sbar}{\ensuremath{\overline{s}}}
\newcommand{\Kn}{\ensuremath{\overline{K_n}^\times}}
\newcommand{\tame}{\operatorname{tame}}
\newcommand{\Qpt}{\ensuremath{\Q_p^{\tame}}}
\newcommand{\Qpu}{\ensuremath{\Q_p^{\unr}}}
\newcommand{\scrT}{\ensuremath{\mathfrak{T}}}
\newcommand{\That}{\ensuremath{\hat{\mathfrak{T}}}}
\newcommand{\Gal}{\operatorname{Gal}}
\newcommand{\Aut}{\operatorname{Aut}}
\newcommand{\tors}{\operatorname{tors}}
\newcommand{\Zhat}{\hat{\Z}}
\newcommand{\linf}{\ensuremath{l_\infty}}
\newcommand{\Lie}{\operatorname{Lie}}
\newcommand{\GL}{\operatorname{GL}}
\newcommand{\End}{\operatorname{End}}
\newcommand{\aone}{\ensuremath{(a_1,\ldots,a_k)}}
\newcommand{\raone}{\ensuremath{r(a_1,\ldots,a_k,N)}}
\newcommand{\rtwoplus}{\ensuremath{\R^{2  +}}}
\newcommand{\rkplus}{\ensuremath{\R^{k +}}}
\newcommand{\length}{\operatorname{length}}
\newcommand{\Vol}{\operatorname{Vol}}
\newcommand{\cross}{\operatorname{cross}}
\newcommand{\GoN}{\Gamma_0(N)}
\newcommand{\GeN}{\Gamma_1(N)}
\newcommand{\GAG}{\Gamma \alpha \Gamma}
\newcommand{\GBG}{\Gamma \beta \Gamma}
\newcommand{\HGD}{H(\Gamma,\Delta)}
\newcommand{\Ga}{\mathbb{G}_a}
\newcommand{\Div}{\operatorname{Div}}
\newcommand{\Divo}{\Div_0}
\newcommand{\Hstar}{\cal{H}^*}
\newcommand{\txon}{\tilde{X}_0(N)}
\newcommand{\sep}{\operatorname{sep}}
\newcommand{\notp}{\not{p}}
\newcommand{\Aonek}{\mathbb{A}^1/k}
\newcommand{\Wa}{W_a/\mathbb{F}_p}
\newcommand{\Spec}{\operatorname{Spec}}

\newcommand{\abcd}{\left[ \begin{array}{cc}
a & b \\
c & d
\end{array} \right]}

\newcommand{\abod}{\left[ \begin{array}{cc}
a & b \\
0 & d
\end{array} \right]}

\newcommand{\unipmatrix}{\left[ \begin{array}{cc}
1 & b \\
0 & 1
\end{array} \right]}

\newcommand{\matrixeoop}{\left[ \begin{array}{cc}
1 & 0 \\
0 & p
\end{array} \right]}

\newcommand{\w}{\omega}
\newcommand{\Qpi}{\ensuremath{\Q(\pi)}}
\newcommand{\Qpin}{\Q(\pi^n)}
\newcommand{\pibar}{\overline{\pi}}
\newcommand{\pbar}{\overline{p}}
\newcommand{\lcm}{\operatorname{lcm}}
\newcommand{\trace}{\operatorname{trace}}
\newcommand{\OKv}{\mathcal{O}_{K_v}}
\newcommand{\Abarv}{\tilde{A}_v}
\newcommand{\Kbar}{\overline{K}}
\newcommand{\pl}{\rho_l}
\newcommand{\plt}{\tilde{\pl}}
\newcommand{\plo}{\pl^0}
\newcommand{\Du}{\underline{D}}
\newcommand{\A}{\mathbb{A}}
\newcommand{\D}{\underline{D}}
\newcommand{\op}{\operatorname{op}}
\newcommand{\Glt}{\tilde{G_l}}
\newcommand{\gl}{\mathfrak{g}_l}
\newcommand{\gltwo}{\mathfrak{gl}_2}
\newcommand{\sltwo}{\mathfrak{sl}_2}
\newcommand{\h}{\mathfrak{h}}
\newcommand{\tA}{\tilde{A}}
\newcommand{\sss}{\operatorname{ss}}
\newcommand{\X}{\Chi}
\newcommand{\ecyc}{\epsilon_{\operatorname{cyc}}}
\newcommand{\hatAl}{\hat{A}[l]}
\newcommand{\sA}{\mathcal{A}}
\newcommand{\sAt}{\overline{\sA}}
\newcommand{\OO}{\mathcal{O}}
\newcommand{\OOB}{\OO_B}
\newcommand{\Flbar}{\overline{\F_l}}
\newcommand{\Vbt}{\widetilde{V_B}}
\newcommand{\XX}{\mathcal{X}}
\newcommand{\GbN}{\Gamma_\bullet(N)}
\newcommand{\Gm}{\mathbb{G}_m}
\newcommand{\Pic}{\operatorname{Pic}}
\newcommand{\FPic}{\textbf{Pic}}
\newcommand{\solv}{\operatorname{solv}}
\newcommand{\Hplus}{\mathcal{H}^+}
\newcommand{\Hminus}{\mathcal{H}^-}
\newcommand{\HH}{\mathcal{H}}
\newcommand{\Tr}{\operatorname{Tr}}
\newcommand{\gK}{\mathfrak{g}_K}
\newcommand{\TT}{\mathcal{T}}
\newcommand{\TTT}{\ilim T_i}
\newcommand{\grad}{\operatorname{grad}}
\newcommand{\Jac}{\operatorname{Jac}}
\newcommand{\loc}{\operatorname{loc}}
\newcommand{\PSL}{\operatorname{PSL}}

\begin{abstract}
Given an algebraic curve $C_{/\Q}$ having points everywhere
locally and endowed with a suitable involution, we show that there
exists a positive density family of prime quadratic twists of $C$
violating the Hasse principle. The result applies in particular to
$w_N$-Atkin-Lehner twists of most modular curves $X_0(N)$ and to
$w_p$-Atkin-Lehner twists of certain Shimura curves $X^{D+}$.
\end{abstract}

\section{Introduction}
\subsection{Some motivation}
Let $C_{/\Q}$ be a (nonsingular, projective, geometrically
integral) algebraic curve.  We say $C$ \emph{violates the Hasse
principle} if for all places $p \leq \infty$, $C(\Q_p) \neq
\emptyset$, but $C(\Q) = \emptyset$.  It is of interest to
\emph{construct} families of such curves in a systematic way, and
also to \emph{discover} naturally occurring families.\footnote{To
be sure, the existence of a clear boundary between construction
and discovery seems dubious, and indeed some of the results
presented here toe this putative line.} Although one has various
results showing that there are ``many'' curves violating the Hasse
principle -- e.g, infinitely many for every genus $g \geq 1$
\cite{CM} -- to the best of my knowledge the literature contains
no ``natural'' infinite\footnote{For an interesting finite family
constructed by a systematic method, see \cite{RSY}.} family of
curves $C_{/\Q}$ provably violating the Hasse principle, by which
I mean an infinite sequence of curves of prior
arithmetic-geometric interest and not just constructed for this
purpose.\footnote{The situation is different if one allows the
possibility of a variable base extension: see
\cite{HasseShimura}.}  Our primary goal in this note is to exhibit
such a ``natural'' infinite family.
\subsection{Statements of the main results} \textbf{} \\
Let $N$ be a squarefree positive integer, $N \neq p \equiv 1 \pmod
4$ a prime, and let $C(N,p)_{/\Q}$ be the \emph{quadratic twist}
of $X_0(N)$ by the Atkin-Lehner involution $w_N$ and the quadratic
extension $\Q(\sqrt{p})/\Q$.  (Precisely what this means will be
reviewed shortly.) These curves are moduli spaces of elliptic
$\Q$-curves (e.g. \cite{Ellenberg}): roughly speaking,
$C(N,p)(\Q)$ parameterizes elliptic curves defined over
$\Q(\sqrt{p})$ which are cyclically $N$-isogenous to their Galois
conjugates.
\begin{thm}
\label{THM1} For $131 < N \neq 163$ a squarefree integer, the set
of primes $p \equiv 1 \pmod 4$ such that $C(N,p)$ violates the
Hasse principle over $\Q$ has positive density.
\end{thm}
\noindent Upon examining the proof of Theorem \ref{THM1}, it
swiftly became clear that the properties of $X_0(N)$ and $w_N$
needed for the argument could be axiomatized to yield a general
criterion for prime quadratic twists violating the Hasse
principle.  This generalization costs nothing extra -- indeed, it
seems if anything to clarify matters -- and it shall turn out to
have (at least) one other interesting application.
\\ \\
The axiomatic version goes as follows: let $C_{/\Q}$ be a curve
and $\iota: C \ra C$ be a $\Q$-rational involution.  Let
$\Q(\sqrt{d})/\Q$ be a quadratic extension, with nontrivial
automorphism $\sigma_d$.  Then there is a curve $C_d =
\mathcal{T}(C,\iota,\Q(\sqrt{d})/\Q)$ called the \emph{quadratic
twist} of $C$ by $\iota$ and $\Q(\sqrt{d})/\Q$ (or just by $d$).
The curve $C_d$ is isomorphic to $C$ over $\Q(\sqrt{d})$ but not,
in general, over $\Q$: $\sigma_d$ acts on $C_d(\Q(\sqrt{d}))$ by
$P \mapsto \iota(\sigma_d(P))$. Its existence follows from the
principle of Galois descent: the $\overline{\Q}/\Q$-twisted forms
of $C$ with respect to the automorphism group generated by $\iota$
are parameterized by
\[H^1(\overline{\Q}/\Q,\langle \sigma \rangle) =
\Hom(\Gal(\overline{\Q}/\Q), \pm 1) \cong \Q^{\times}/\Q^{\times
2}.\]
\begin{thm}
\label{THM2} Let $C_{/\Q}$ be an algebraic curve.  Suppose there
is a $\Q$-rational
involution $\iota$ on $C$ such that: \\
(i) $\{P \in C(\Q) \ | \ \iota(P) = P\} = \emptyset$. \\
(ii) There exists $P_0 \in C(\overline{\Q})$ such that $\iota(P_0)
= P_0$. \\
(iii) For all $\ell \leq \infty$, $C(\Q_{\ell}) \neq \emptyset$,
i.e., $C$ has
points everywhere locally. \\
(iv) The quotient $C/\iota$ has finitely many $\Q$-rational
points. \\
Then the set of primes $p \equiv 1 \pmod 4$ for which the twist
$C_p = \mathcal{T}(C,\iota,\Q(\sqrt{p})/\Q)$ violates the Hasse
principle has positive density.
\end{thm}
\noindent Remark 1.1: Each of the assumptions (i) and (ii) is
necessary for the conclusion to hold: if (i) fails then every
quadratic twist $C_d$ has $\Q$-points; whereas if (ii) fails then
the set of squarefree $d$ for which $C_d$ has points everywhere
locally is finite \cite[Prop. 5.3.2]{Skorobogatov}. One must
assume at least (iii$'$): $C_{/\Q_{\ell}}$ has quadratic points
for all $\ell \leq \infty$ (or, in the terminology of
\cite{HasseShimura}, that $m_{\loc}(C) \leq 2$).  It is surely not
the case that (iii$'$) is \emph{in general} sufficient for the
existence of twists $C_d$  (prime or otherwise) violating the
Hasse principle, but it might be interesting to try to modify the
argument so as to apply to some particular curves satisfying
(iii$'$), e.g. certain Shimura curves $X^D_0(N)$. Finally, some
additional hypothesis, like (iv), is needed to ensure that $C$ is
not the projective line!
\\ \\
Remark 1.2: A pair $(C,\iota)$ can satisfy the hypotheses of
Theorem \ref{THM2} only if $C/\iota$ has genus at least one and
$C$ has genus at least two.  Indeed, since $C$ has points
everywhere locally, so does $C/\iota$, so if it had genus zero it
would -- by the Hasse principle -- be isomorphic to $\PP^1$ and
have infinitely many $\Q$-points.  But $C/\iota$ would have genus
zero if $C$ had genus zero (clearly) or genus one (since $\iota$
has fixed points).
\\ \\
Remark 1.3: From the proof one can deduce an explicit lower bound
on the density of the set $\mathcal{P} = \mathcal{P}(C,\iota)$ of
primes such that $C_p$ violates the Hasse principle in terms of
the genera of $C$ and $C/\iota$.  However, it does \emph{not} give
us an effective way to find any elements of $\mathcal{P}$ unless
we can find all the $\Q$-points on $C/\iota$.  In the case of
$X^+_0(N) = X_0(N)/w_N$, this is a notoriously difficult open
problem, c.f. $\S 3.3$.
\\ \\
One may ask why in (iii) we do not just assume that $C$ has a
$\Q$-point: otherwise our curve \emph{already} violates the Hasse
principle, and twisting it to get further violations seems less
interesting. The advantage of stating the result the way we have
is that it applies \emph{even if we don't know whether or not} $C$
has $\Q$-rational points.  This distinction is illustrated in the
following additional application of Theorem \ref{THM2}.
 \\ \\
Let $D = p_1 \cdots p_{2r}$ be a nontrivial squarefree product of
an even number of primes, and consider the Shimura curve
$X^D_{/\Q}$ (e.g. \cite[Chapter 0]{Thesis}).  Unfortunately (iii)
does not hold -- e.g. $X^D(\R) = \emptyset$
\cite{Ogg2}\footnote{The result is originally due to Shimura.} --
so Theorem 2 does not apply. On the other hand, the full
Atkin-Lehner group $W$ consists of $2^{2r}$ commuting involutions
$w_d$, one for each $1 \leq d \ | \ D$.  It turns out that $X^{D+}
:= X^D/w_D$ has points everywhere locally \cite[Main Theorem
2]{Thesis}. Since $W$ is commutative of order at least $4$, for a
prime $q \ | \ D$, $w_q$ induces a well-defined, nontrivial
involution on $X^{D+}$, which we continue to denote by $w_q$. Let
$C_p(D,w_q)$ be the twist of $X^{D+}$ by $w_q$ and
$\Q(\sqrt{p})/\Q$.
\begin{thm}
\label{THM3} There exists\footnote{One can certainly compute an
explicit $D_0$; for lack of any application I have not done so
here.} an integer $D_0$ with the following property: for pairwise
distinct primes $q
> 163, \ p_2, \ldots,\  p_{2r}$ such that: \\
(a)  $(\frac{q}{p_i}) \neq 1$ for all $i$; and \\
(b) $D := q \cdot p_2 \cdots p_{2r} > D_0$, \\ the set of primes
$p$ for which $C_p(D,w_q)$ violates the Hasse principle over $\Q$
has positive density.
\end{thm}
\noindent Although this family is perhaps less ``natural'' than
that of Theorem \ref{THM1} -- the modular interpretation of
$C_p(D,w_q)$ is rather abstruse -- it is interesting for other
reasons, as we explain at the end.
\\ \\
The proofs are given in $\S 2$. In $\S 3$ we discuss certain
complements, in particular a generalization of Theorem \ref{THM2}
to twists by automorphisms of prime degree $p$ on curves defined
over a number field containing the $p$th roots of unity.
\subsection{Connections with the Inverse Galois Problem for $PSL_2(\F_p)$}
Our motivation for considering Atkin-Lehner twists of $X_0(N)$ by
\emph{primes} $p$ in particular comes from work of K.-y. Shih, who
showed that if $(\frac{N}{p}) = -1$, then there is a covering of
curves $Y \ra C(N,p)$, Galois over $\Q$ with group $PSL_2(\F_p)$
\cite[Thm. 8]{Shih1}.  There is thus an ulterior motive for
studying the locus $C(N,p)(\Q)$: each such point $P$ yields by
specialization a homomorphism $\rho_P: \Gal(\overline{\Q}/\Q) \ra
PSL_2(\F_p)$. So if $\rho_P$ is surjective for some $P$, we get a
realization of $PSL_2(\F_p)$ as a Galois group over $\Q$ (an open
problem, in general).
 \\ \indent If we have such a surjective specialization, let us say that ``Shih's strategy succeeds'' for
these values of $N$ and $p$.  Remarkably -- more than $30$ years
after \cite{Shih1} -- Shih's strategy remains the state of the
art: with the single exception of a theorem of Malle \cite{Malle}
dealing with the case $(\frac{5}{p}) = -1$ -- every
realization of $PSL_2(\F_p)$ as a Galois group over $\Q$ is an instance of the success of Shih's strategy.  \\
\indent Shih himself gave a complete analysis of the cases where
$C(N,p)$ has genus zero. In \cite{GG}, following up on a
suggestion of Serre, I analyzed some of the cases of genus one. At
the end of that paper I asked three questions about rational
points on $C(N,p)$, the last being whether there could exist
points everywhere locally but not globally.  In an early draft of
this paper I had billed Theorem \ref{THM1} as an answer to this
question, but it is not, really, since the primes constructed in
the proof all satisfy $(\frac{N}{p}) = 1$. So I have included an
Appendix to this paper in which the topic of local and global
points on $C(N,p)$ under the hypothesis $(\frac{N}{p}) = -1$ is
revisited. We make some further remarks and questions about the
genus one cases, and, especially, we modify the argument of
Theorem \ref{THM2} to give many Hasse principle violations when
$(\frac{N}{p}) = -1$, showing that there are purely global
obstructions to the success of Shih's method.

\section{Proofs} \noindent \emph{Proof} of Theorem 2: For a squarefree $d \neq 1$,
there are natural set maps  \[\alpha_d: C_d(\Q) \hookrightarrow
C(\Q(\sqrt{d}))\] and \[\beta_d: C(\Q(\sqrt{d})) \ra
(C/\iota)(\Q(\sqrt{d})). \] Put
\[S_d = (\beta_d \circ \alpha_d)(C_d(\Q)). \] Then $S_d \subset
(C/\iota)(\Q)$. Moreover, $(C/\iota)(\Q) = \bigcup_d S_d \cup
\iota(C(\Q))$, and for $d \neq d'$, $P \in S_d \cap S_{d'}$
implies that $P \in C(\Q(\sqrt{d})) \cap C(\Q(\sqrt{d'})) =
C(\Q)$.  But $S_d \cap C(\Q)$ consists of $\Q$-rational
$\iota$-fixed points, which we have assumed in (i) do not exist,
so that for $d \neq d'$, $S_d \cap S_{d'} = \emptyset$. By (iv),
$(C/\iota)(\Q)$ is finite, and we conclude that the set of $d$ for
which $S_d \neq \emptyset$ is finite.  In particular, for all
sufficiently large primes $p$, the prime twists $C_p$ have no
$\Q$-rational points. \\ \indent Thus it shall suffice to
construct a set of primes $p \equiv 1 \pmod 4$ of positive density
such that $C_p$ has $\Q_{\ell}$-rational points for all $\ell \leq
\infty$. A key observation is that if $p$ is a square in
$\Q_{\ell}$, then since $\Q_{\ell}$ contains $\Q(\sqrt{p})$,
$C_p(\Q_{\ell}) = C(\Q_{\ell})$, which is nonempty by (iii). In
particular, since $p
> 0$, $C_p(\R) \neq \emptyset$ for all $p$.
\\ \indent Let $M_1$ be a positive integer such that for $\ell
> M_1$, $C$ extends to a smooth relative curve
$C_{/\Z_{\ell}}$. If $\ell > M_1$ is prime to $p$, then we claim
that $C_p$ also extends smoothly over $\Z_{\ell}$.  Indeed, the
extension $\Q_{\ell}(\sqrt{p})/\Q_{\ell}$ is unramified (note that
we use $p \equiv 1 \pmod 4$ here), and after this base change
$C_p$ becomes isomorphic to $C$.  As above, (iii) and (iv) imply
that $C$ (hence also $C_p$) has positive genus, so $C$ admits a
minimal regular $\Z_{\ell}$-model. But it is known that formation
of the
minimal regular model commutes with unramified base change, so the minimal model of $C_{/\Z_{\ell}}$ must already have been smooth. \\
\indent Let $g$ be the genus of $C$; notice that it is also the
genus of $C_p$ for all $p$.  By the Weil bounds for curves over
finite fields, there exists a number $M_2$ such that if $\ell >
M_2$, every nonsingular curve $C_{/\F_{\ell}}$ of genus $g$ has an
$\F_{\ell}$-rational point.  Thus, if $\ell > M = \max\{M_1,M_2\}$
and is different from $p$, then $C_p$ admits a regular $\Z_{\ell}$
model whose special fiber has a smooth $\F_{\ell}$-rational point;
by Hensel's Lemma this implies that $C_p(\Q_{\ell}) \neq
\emptyset$. \\ \indent It remains to choose $p$ to take care of
the primes $\ell \leq M$ and $\ell = p$.  In the former case we
may just assume that $p \equiv 1 \pmod 8$ and that $p$ is a
quadratic residue modulo every odd $\ell \leq M$, so that by the
above observation we get that $C_p(\Q_{\ell}) = C(\Q_{\ell}) \neq
\emptyset$.  Finally, to get that $C_p(\Q_p) \neq \emptyset$, we
use (ii) the existence of a geometric $\iota$-fixed point $P_0$.
If we choose $p$ to split completely in $\Q(P_0)$, then $P_0$ is a
$\Q_p$-rational fixed point of $\iota$, so is an element of
$C(\Q_p) \cap C_p(\Q_p)$. In all we have required $p$ to split
completely in a finite collection of number fields, so all primes
splitting completely in the compositum will do. Cebotarev's
theorem implies that this set of primes has positive density.
\\ \\
\emph{Proof} of Theorem \ref{THM1}: We just check that the
hypotheses of Theorem \ref{THM2} apply with $C = X_0(N)$ and
$\iota = w_N$. It is well-known that there are always $w_N$-fixed
points -- so (i) holds -- and that
\[\min \ \{ [\Q(P):\Q] \ | \ P \in X_0(N)(\overline{\Q}), \ w_N(P) =
P\} =  h(\Q(\sqrt{-N})), \] the class number of $\Q(\sqrt{-N})$
(e.g. \cite[Prop. 3]{Ogg}).  (We use here that $N$ is squarefree.)
Thus (ii) holds unless $h(-N) = 1$, i.e., unless $N = 1, \ 2, \ 3,
\ 7,\ 11,\ 19,\ 43, \ 67$ or $163$. Using the genus formula for
$X^+_0(N)$ (and the standard upper bound on the class numbers of
imaginary quadratic fields), it is straightforward to compute the
complete list of $N$ for which $X^+_0(N)$ has genus at most one.
We shall not give this list here, but the largest such $N$ is
$131$ (compare with e.g. \cite{Bars}), verifying (iii). Finally,
the rationality of the cusps gives $X_0(N)(\Q) \neq \emptyset$ for
all $N$, hence (iv) holds.
\\ \\
\emph{Proof} of Theorem \ref{THM3}: Again we will verify the
hypotheses of Theorem \ref{THM2}, now with $C = X^{D+}$ and
$\iota$ the image of $w_q$.  The congruence conditions in the
statement of theorem ensure that $w_q$ has fixed points on $X^D$
(e.g. \cite{Ogg2} or \cite[Prop. 48]{Thesis}); \emph{a fortiori}
its image on $X^{D+}$ has fixed points, so (i) holds. The field of
definition $\Q(P)$ of any $w_q$-fixed point $P$ contains the
Hilbert class field of $\Q(\sqrt{-q})$, so if $q > 163$,
$[\Q(P):\Q] > 2$.  The degree of the field of definition of the
image of $P$ on the involutory quotient $X^{D+}$ is at least
$\frac{1}{2}[\Q(P):\Q]$, so there are no $\Q$-rational
$\iota$-fixed points on $X^{D+}$.\footnote{Note here the analogy
between
$X_0(N)$ and $X^{D+}$, rather than $X^D$.} Thus (ii) holds. \\
\indent For (iv), it is enough to know that the genus of
$X^D/\langle w_D, w_q \rangle$ is at least $2$ for all
sufficiently large $D$.  But a routine calculation using the
formulae for the genus of $X^D$ and the number of fixed points of
the $w_d$'s shows that even the genus of the full Atkin-Lehner
quotient $X^D/W$ approaches $\infty$ with $D$ (e.g.
\cite[Corollary 50]{Thesis}). \\ \indent As already mentioned, I
showed that (iii) holds for \emph{all} squarefree $D$ in my
(unpublished) Harvard thesis \cite[Main Theorem 2]{Thesis}. In the
meantime, Rotger, Skorobogatov and Yafaev have proved a more
general result \cite[Theorem 3.1]{RSY}. Both proofs use a result
of Ogg for $\ell = \infty$; a trace formula for $(\ell,D) = 1$;
and the Cerednik-Drinfeld uniformization for $\ell \ | \ D$; i.e.,
they are similar enough so that I feel no need to reproduce the
details of my argument here.
\\ \\
\section{Complements}

\subsection{Variants of Theorem \ref{THM2}} \textbf{} \\
Theorem \ref{THM2} extends immediately to the case of $(C,\iota)$
defined over an arbitrary number field $K$, still twisting by the
quadratic (for all but finitely many $p$) extensions
$K(\sqrt{p})/K$.  Of course, these need not be ``prime'' quadratic
twists with respect to $K$, but this can be remedied.
\\ \\
Indeed, there is an analogue of Theorem \ref{THM2} for a curve $C$
over a number field $K \supset \Q(\mu_p)$ endowed with a
$K$-rational automorphism $\varphi$ of prime order $p$.  Kummer
theory gives:
\[K^{\times}/K^{\times p} = H^1(K,\langle \varphi \rangle). \] By
Galois descent, an element $x \in K^{\times}$ gives rise to a
twist $C_x = \mathcal{T}(C,\varphi,K(x^{\frac{1}{p}})/K)$.
\newcommand{\Fix}{\operatorname{Fix}}
\begin{thm}
Let $K \supset \Q(\mu_p)$ be a number field, and let $C_{/K}$ be
an algebraic curve. Suppose there is a $K$-rational
automorphism $\varphi: C \ra C$ of prime order $p$ such that: \\
(i) $\{P \in C(K) \ | \ \varphi(P) = P\} = \emptyset$. \\
(ii) There exists $P_0 \in C(\overline{K})$ such that
$\varphi(P_0) = P_0$. \\
(iii) $C$ has points everywhere locally. \\
(iv) The quotient $C/\varphi$ has finitely many $K$-rational
points. \\
Then there exists a modulus $\mathfrak{m}$ of $K$ such that the
prime ideals $\mathfrak{p}$ of $\mathfrak{o}_K$ which are
generated by an element $\pi \equiv 1 \pmod{\mathfrak{m}}$ such
that $C_{\pi}$ violates the Hasse principle over $K$ have positive
density.
\end{thm}
\noindent Here $\mathfrak{m}$ is chosen so that $\pi \equiv 1
\pmod{\mathfrak{m}}$ implies that $\pi$ is totally positive and a
perfect $p$-th power in the completion of $K$ at any prime
$\mathfrak{p}$ lying over $p$; e.g., when $K = \Q$ we take
$\mathfrak{m} = 8 \cdot  \infty$.  The remainder of the proof is
left to the reader.
\subsection{Rational points on Atkin-Lehner quotients of
Shimura curves} Let us end by calling attention to the sequence of
curves $X^{D+}$: for me it is the example \emph{par excellence} of
a naturally occurring family of curves with points everywhere
locally.  In fact work of Jordan, Rotger and others shows that the
modular interpretation of $X^{D+}$ is in many respects more
natural than that of $X^D$.  Recall the ``folk conjecture'' that
for sufficiently large squarefree $N$, $X^+_0(N)(\Q)$ consists
only of cusps and CM points.
  As alluded to above, this is an extremely difficult conjecture,
since the favorable case for determining the $\Q$-points on
$C_{/\Q}$ is when $\Jac(C)$ has a $\Q$-factor of rank less than
its dimension. But analysis of the sign in the functional
equations shows -- assuming the conjecture of Birch and
Swinnerton-Dyer (henceforth BSD) -- that this \emph{never} happens
for $X^+_0(N)$.
\\ \indent By work of Jacquet-Langlands and Faltings,
the same holds for $X^{D+}$, so that the study of $X^{D+}(\Q)$ is
again very difficult.  The $X^{D+}$ version of the above folk
conjecture is that for sufficiently large $D$, $X^{D+}(\Q)$
consists only of CM points. But the CM points are well understood,
and one can see in particular that the set of $D$ such that
$X^{D+}$ has a $\Q$-rational CM point has density
zero.\footnote{This uses the zero density of the set of integers
with a bounded number of prime divisors; among discriminants $D =
p_1 p_2$, the probability of rational CM point is
$(1-(1-\frac{1}{4})^9) \approx .925$.}  Thus either (I) the curves
$X^{D+}$ violate the Hasse principle for a density one set of
squarefree integers $D$; (II) there is some (as yet unknown)
specific phenomenon which puts many $\Q$-rational points on
Atkin-Lehner quotients of Shimura curves; or (III) our
conventional wisdom about $\Q$-points on algebraic curves -- i.e.,
that they are ``typically'' relatively sparse unless there is some
good reason -- is completely wrong. Each of three options is
fascinating in its own way, but which is true?!?
\\ \indent
It seems to me that Theorems \ref{THM2} and \ref{THM3} provide
some evidence against (III).  But the curves $X^{D+}$ are hardly
``randomly chosen'': the conjectured nonexistence of non-CM
$\Q$-rational points implies the nonexistence of certain kinds of
$GL_2$-type abelian surfaces $A_{/\Q}$ and (\emph{a fortiori},
assuming Serre's conjecture) certain kinds of modular forms.  Our
current understanding of these associated objects is so limited
that we certainly cannot dismiss the possibility of (II).
\\ \\
\section*{Appendix: Atkin-Lehner Twists With $(\frac{N}{p}) = -1$}

\noindent This is essentially an addendum to \cite{GG}.  Although
some of the results of \emph{loc. cit.} will be revisited here
with a slightly different emphasis, for more complete accounts the
reader should consult \cite{GG} as well as \cite{Shih1} and
\cite{Serre}.
\\ \\
We will suppose throughout that $N > 1$ is a squarefree integer
and $p$ is an odd prime such that $(\frac{N}{p}) = -1$; by
\cite[Thm. 8]{Shih1} there is then a $\Q$-rational
$PSL_2(\F_p)$-Galois cover $Y \ra C(N,p)$.\footnote{In the body of
the text we defined $C(N,p)$ only for $p \equiv 1 \pmod 4$; indeed
the proof of Theorem \ref{THM2} required $p \equiv 1 \pmod 8$. But
Shih's theorem holds for $p \equiv -1 \pmod 4$ if we define
$C(N,p)$ as the twist of $X_0(N)$ by $w_N$ and $-p$ -- i.e., in
general as the twist by $w_N$ and $p^* = p^{\frac{p-1}{2}}$.}  We
would like to know for which values of $N$ and $p$ there exists a
point $P \in C(N,p)(\Q)$ whose associated homomorphism $\rho_P:
\Gal(\overline{\Q}/\Q) \ra \PSL_2(\F_p)$ is surjective.  In
particular we would like to know whether $C(N,p)(\Q) \neq
\emptyset$, and especially whether there are any non-CM points
(since the specialization homomorphism at a CM point will have
abelian image upon restriction to the corresponding imaginary
quadratic field, it will not be surjective for any odd $p$).
Intuitively, one might expect $\rho_P$ to be surjective for
``most'' non-CM points.  Indeed, work of Hilbert, Faltings and
Serre shows that whenever $C(N,p)(\Q)$ is infinite, there are
infinitely many irreducible specializations.
\\ \\
The cases in which $C(N,p)$ has genus zero ($N = 2$, $3$, $5$,
$7$, $10$, $13$) are decisively treated in \cite{Shih1},
\cite{Serre} and \cite{GG}: in particular, when $N = 2$, $3$ or
$7$ (class number $1$!), Shih's strategy succeeds for all $p$
(with $(\frac{N}{p}) = -1$).  The genus one cases are $N = 11$,
$14$, $15$, $17$, $19$, $21$; for $N = 11$ and $19$ (class number
$1$!) we will say only that Shih's strategy works when $p \equiv 1
\pmod 4$ conditionally on BSD and refer the reader to \cite{GG}
for more details.  To deal with some of the other cases, the
following result was used \cite[Theorem 11]{GG}.\footnote{Let us
note, as we did in \cite{GG}, that the implication
$\Longleftarrow$ of Theorem \ref{DRTHM} is a special case of a
theorem of Gonzalez \cite[Thm 6.2]{Quer}.}
\begin{thm}
\label{DRTHM}
For prime $N$, $C(N,p)(\Q_N) = \emptyset \iff N
\equiv 1 \pmod 4$.
\end{thm}
\noindent   This shows the failure of Shih's method when $N = 17$
(and also $N = 5, \ 13$; a similar argument shows that
$C(10,p)(\Q_5) \neq \emptyset$ iff $(\frac{5}{p}) = 1$).  The
other cases ($N = 14$, $15$, $21$) were not analyzed in \cite{GG}
since they could not lead to any new Galois groups (nor even to
the recovery of Malle's result).
\\ \\
Note the remarkable fact that in every case above, we either had
an obvious $\Q$-rational point or a local obstruction at a prime
$\ell \ | \ N$; in particular there were no violations of the
Hasse principle.
\\ \\
Here we want to point out that the cases $N = 14$ and $N = 21$
nevertheless give rise to some interesting Diophantine problems.
The analysis of local points on $C(14,p)$ and $C(21,p)$ is
complete except at the prime $p$: interestingly, computations
suggest that there are points everywhere locally iff there are
points at every place except $p$, which leads me to believe that
there is a criterion for the emptiness of $C(N,p)(\Q_p)$ which is
equally simple as that of $\Q_{\ell}$ for $\ell \ | N$ and
``quaternionically linked to it.''\footnote{I have, in part,
refrained from serious contemplation of the minimal model of
$C(N,p)_{/\Z_p}$ for hope of it serving as part of some future
thesis project.  So although you are of course welcome to work on
this problem, please let me know if you solve it.} When (for
instance) $N = 14$ and $p \equiv 17 \pmod {56}$, computations
suggest that there are always points everywhere locally. Despite
the fact that there are no ``obvious'' $\Q$-rational points here,
for the first $101$ such primes in this congruence class we do in
fact get an elliptic curve, necessarily of odd analytic rank, so
this gives $101$ cases of the success of Shih's strategy.  It
would be surprising if this phenomenon persisted for all $p$ in
this congruence class. However, our lack of counterexamples is the
``Selmer dual'' of a phenomenon encountered in \cite{GG}: let
$J(14,p)$ be the Jacobian of $C(14,p)$, i.e., the quadratic twist
of $X_0(14)$ by $p^*$ in the usual sense. Since $J_0(14)[2] \cong
X_0(14)[2] \cong \Z/2\Z$, and twists by primes $p \equiv 1 \pmod
4$ have odd analytic rank, then (assuming BSD) the curve $C(14,p)$
can only represent a nontrivial element of $\Sha(\Q,J(14,p))[2]$
if the $2$-Selmer rank of $J(14,p)$ is at least $3$ (and in fact
at least $4$ because the contribution of $\Sha[2]$ to the
$2$-Selmer rank will be even). Now recall that in \cite{GG} we
were unable to find prime twists of $X_0(11)$ or $X_0(19)$ of rank
at least $3$!  It would be very interesting to have some
conjectural asymptotics on the variation of $2$-Selmer ranks in
families of prime twists that would give us a hint as to how far
we ought to look before being surprised by the lack of examples in
each case.
\\ \\
I am not aware of a single example of the success of Shih's
strategy when $C(N,p)$ has genus at least two.  These would
necessarily come from exceptional $\Q$-rational points on
$X^+_0(N)$ with the additional condition that the quadratic field
of the preimage is $\Q(\sqrt{p^*})$ for some prime $p$ with
$(\frac{N}{p}) = -1$ (for the associated $\rho_P$ must be
surjective, but this seems to be the least of our
worries).\footnote{I have not, however, had the opportunity to
read through all of the rapidly growing literature on $\Q$-curves
in search of such examples.  It would be very useful if there
existed an online database containing all known exceptional
rational points on $X^+_0(N).$} Note that the ``folk conjecture''
of $\S 3.2$ predicts that Shih's strategy succeeds for only
finitely many pairs $(N,p)$ with $g(C(N,p)) \geq 2$, so in
particular that we will not be able to obtain all $PSL_2(\F_p)$'s
as a Galois group over $\Q$ by this method.
\\ \\
We shall now show that there are plenty of curves $C(N,p)$ with $p
\equiv 1 \pmod 4$, $(\frac{N}{p}) = -1$ violating the Hasse
principle.  Assume for simplicity that $N$ is \emph{prime}. Assume
also that $N>163$, so that $C(N,p)(\Q) = \emptyset$ for all but
finitely many primes $p$.  We will try to run through the argument
of Theorem $2$ except with the congruence condition $(\frac{N}{p})
= -1$, i.e., $p$ is inert in $\Q(\sqrt{N})$.  There are however
two issues to be addressed: the first is that $C(N,p)$ need not
have $\Q_N$-rational points; indeed by Theorem \ref{DRTHM} we know
that it will iff $N \equiv -1 \pmod 4$, so let us assume this
condition on $N$.  It is then indeed the case that for all but
finitely many $p$ satisfying the modified conditions, $C_p$
violates the Hasse principle. However, there is a second issue:
since one of our conditions on $p$ is not a splitting condition,
it is not \emph{a priori} clear that our conditions on $p$ are
consistent, i.e., correspond to some nonvoid Cebotarev set.
\\ \indent Let us now carefully check the consistency.   To be
precise, we are imposing the following conditions on $p$: (a) $p$
splits completely in $K_1 := \Q(\zeta_8)$ -- i.e., $p \equiv 1
\pmod 8$; (b) $p$ is inert in $\Q(\sqrt{N})$; (c) $p$ splits in
sufficiently many quadratic fields $\Q(\sqrt{\ell^*})$ for $\ell$
odd and different from $p$ and $N$ (namely, such that for the
remaining primes $\ell$, $C(N,p)$ has smooth reduction modulo
$\ell$ and an $\F_{\ell}$-rational point on its special fiber);
and (d) $p$ splits completely in $\Q(P_0)$, where $P_0$ is a
$w_N$-fixed point corresponding to the \emph{maximal} order of
$\Q(\sqrt{-N})$ (and not the one of conductor $2$).  Let $K_2$ be
the compositum of $\Q(\sqrt{\ell^*})$ for the finite set of odd
primes $\ell$ considered above, and let $K_3 := \Q(P_0,
\sqrt{-N})$, the Hilbert class field of $\Q(\sqrt{-N})$. The
fields $K_i$ for $1 \leq i \leq 3$ are each Galois over $\Q$ and
are mutually linearly disjoint, since their sets of ramified
finite primes are pairwise disjoint. Put $K = K_1 \cdot K_2 \cdot
K_3$, so that $\Gal(K/\Q) = \prod_{i=1}^3 \Gal(K_i/\Q)$.  In $G_3$
there is a unique conjugacy class $C$ consisting of elements
$\sigma$ with nontrivial restriction to $\Q(\sqrt{-N})$.  So, by
Cebotarev, the primes $p$ which are unramified in $K$ and with
corresponding Frobenius automorphism lying in the conjugacy class
$(1,1,C)$ have positive density.  We have shown:
\begin{thm}
Suppose $N > 163$ is prime and congruent to $-1 \pmod 4$.  The set
of primes $p \equiv 1 \pmod 4$, $(\frac{N}{p}) = -1$ such that
$C(N,p)$ violates the Hasse principle over $\Q$ has positive
density.
\end{thm}
\noindent Thus we get an affirmative answer to the question asked
at the end of \cite{GG}: there are ``truly global'' obtructions to
the success of Shih's method.
\\ \\
Acknowledgements: Thanks to Jordan Ellenberg and Bjorn Poonen for
useful comments.  I am grateful to the authors of \cite{RSY} for
their citation of my unpublished thesis work.  The title of this
paper was chosen, in part, to please Siman Wong.
\\ \\
{ \bf Postscript: } Recent (and recently remembered)
correspondences with Dr. Nick Rogers and Noam Elkies suggest that
the computationally observed phenomenon of $2$-Selmer rank at most
$2$ in \emph{prime} twist families of the elliptic curves $X_0(N)$
for $N = 11, \ 19$ (in \cite{GG}) and $N = 14$ (here) may have a
relatively simple explanation -- and in particular, may persist
for all primes in the given congruence classes -- by means of a
complete $2$-descent.  This suggests the possibility -- quite
surprising when compared to Theorem \ref{THM1} -- that when
$X_0(N)$ has genus one there are \emph{no} prime $w_N$-quadratic
twists which violate the Hasse principle. These results may find
their way into the ultimate version of this paper, or they may
appear -- together with much of the material of the Appendix --
elsewhere.

\end{document}